\documentclass{amsart}

\usepackage{amssymb}
\usepackage{mathrsfs}

\theoremstyle{plain}
\newtheorem{theorem}{Theorem}[section]
\newtheorem{lemma}[theorem]{Lemma}

\newtheorem{proposal}[theorem]{Proposal}

\theoremstyle{definition}

\theoremstyle{remark}

\numberwithin{equation}{section}

\begin{document}
\newcommand{\Z}{\mathbb{Z}}
\newcommand{\Q}{\mathbb{Q}}
\newcommand{\R}{\mathbb{R}}
\newcommand{\C}{\mathbb{C}}
\newcommand{\divides}{\mid}
\newcommand{\doesnotdivide}{\nmid}
\newcommand{\abs}[1]{\lvert{#1}\rvert}
\newcommand{\partition}[1]{\ensuremath{\left\langle{#1}\right\rangle}}
\newcommand{\floor}[1]{\left\lfloor{#1}\right\rfloor}
\newcommand{\modtwo}[1]{\left\{{#1}\right\}}
\newcommand{\num}{\nu}
\newcommand{\rem}[1]{R({#1})}
\newcommand{\M}{\mathfrak{M}}
\newcommand{\m}{\mathfrak{m}}
\newcommand{\GEQ}{\succcurlyeq}
\newcommand{\LEQ}{\preccurlyeq}

\title{On the class of dominant and subordinate products}

\author{Alexander Berkovich}
\address{Department of Mathematics, University of Florida, Gainesville,
Florida 32611-8105}
\email{alexb@ufl.edu}          

\author{Keith Grizzell}
\address{Department of Mathematics, University of Florida, Gainesville,
FL 32611-8105}
\email{grizzell@ufl.edu}          


\date{\today}  



\begin{abstract}
In this paper we provide proofs of two new theorems that provide a broad class of partition inequalities and that illustrate a na\"ive version of Andrews' anti-telescoping technique quite well. These new theorems also put to rest any notion that including parts of size 1 is somehow necessary in order to have a valid irreducible partition inequality. In addition, we prove (as a lemma to one of the theorems) a rather nontrivial class of rational functions of three variables has entirely nonnegative power series coefficients.
\end{abstract}

\maketitle

\section{Introduction}

When examining two $q$-products $\Pi_1$ and $\Pi_2$ and their corresponding $q$-series (assuming $|q|<1$), it sometimes happens that the coefficients in the $q$-series for $\Pi_1$ are never less than the coefficients in the $q$-series for $\Pi_2$. When that happens, we say that $\Pi_1$ is \emph{dominant} (in this pair of products) and that $\Pi_2$ is \emph{subordinate}, and we express this relationship with the more succinct notation $\Pi_1 \GEQ \Pi_2$. (Note that $\GEQ$ yields a partial ordering on the set of $q$-products if we identify products that produce the same $q$-series; then, any given product may be dominant when paired with some products, subordinate when paired with others, neither when paired with still other products, and both dominant and subordinate only when paired with ``itself''.) Immediately from this definition it follows that if $\Pi_1 \GEQ \Pi_2$, then the $q$-series determined by $\Pi_1-\Pi_2$ must have nonnegative coefficients, i.e.\ $\Pi_1-\Pi_2 \GEQ 0$. Thus, determining whether a 
given pair of products is a dominant/subordinate pair solves an equivalent positivity problem.

Using the standard notations \cite{An1}
\begin{align*}
(a;q)_L &= \begin{cases}
              1 & \text{if } L=0,\\
              \prod_{j=0}^{L-1}(1-aq^j) &\text{if } L>0,
           \end{cases}\\
(a_1,a_2,\dots,a_n;q)_L &= (a_1;q)_L (a_2;q)_L \cdots (a_n;q)_L,\\
\intertext{and}
(a;q)_\infty &= \lim_{L\to\infty} (a;q)_L,
\end{align*}
we may say that, for example, in the Rogers-Ramanujan difference
\begin{equation}\label{RRdiff}
\frac{1}{(q,q^4;q^5)_{\infty}}- 
\frac{1}{(q^2,q^3;q^5)_{\infty}} \GEQ 0
\end{equation}
the first product is dominant and the second product is subordinate. 
At the 1987 A.M.S.\ Institute on Theta Functions, Ehrenpreis asked if one can prove this dominance without resorting to the Rogers-Ramanujan identities.
In 1999, Kadell \cite{KK1} provided an affirmative answer to this question. 
In 2005, Berkovich and Garvan \cite{BGa} proved a class of finite versions of such inequalities (from which the infinite versions are easily recovered), namely that
\begin{equation}\label{BGadiff}
\frac{1}{(q,q^{m-1};q^m)_L}\GEQ
\frac{1}{(q^r,q^{m-r};q^m)_L}
\end{equation}
if and only if $r \doesnotdivide (m-r)$ and $(m-r) \doesnotdivide r$. Note that this last inequality provides the finite version of \eqref{RRdiff}:
\begin{equation}\label{finiteRRdiff}
\frac{1}{(q,q^4;q^5)_L} \GEQ 
\frac{1}{(q^2,q^3;q^5)_L}.
\end{equation}
In 2011, Andrews \cite{An2} proved the finite little G\"ollnitz inequality
\begin{equation}\label{Andiff}
\frac{1}{(q,q^5,q^6;q^8)_L}\GEQ
\frac{1}{(q^2,q^3,q^7 ;q^8)_L},
\end{equation}
which (in 2012) Berkovich and Grizzell \cite{BGr} generalized to 
\begin{equation}\label{BGrdiff}
\frac{1}{(q,q^{y+2},q^{2y};q^{2y+2})_L}\GEQ
\frac{1}{(q^2,q^{y},q^{2y+1} ;q^{2y+2})_L},
\end{equation}
where $y$ is any odd integer greater than 1.

For \eqref{RRdiff}, \eqref{BGadiff}, and \eqref{BGrdiff}, the proofs in each case relied solely on the construction of a suitable injection. For \eqref{Andiff}, however, Andrews relied primarily on his anti-telescoping technique. A na\"ive version of Andrews' anti-telescoping technique begins with two sequences of products, $\{P(i)\}_{i=1}^\infty$ and $\{Q(i)\}_{i=1}^\infty$, and the desire to show that, for every $L\geq1$,
\begin{equation*}
\frac{1}{P(L)} \GEQ \frac{1}{Q(L)}.
\end{equation*}
One then simply writes
\begin{align}
\frac{1}{P(L)} - \frac{1}{Q(L)} &= \sum_{i=1}^L \frac{Q(i-1)}{P(i)Q(L)}\left(\frac{Q(i)}{Q(i-1)}-\frac{P(i)}{P(i-1)}\right)\label{antitelxy}\\
&= \sum_{i=1}^L \frac{\frac{Q(i)}{Q(i-1)}-\frac{P(i)}{P(i-1)}}{P(i)\cdot\frac{Q(L)}{Q(i-1)}},\label{antitelxyz}
\end{align}
and if one is lucky enough that each addend in \eqref{antitelxyz} is $\GEQ 0$, then that is all one needs to show in order to prove the desired inequality. This bit of serendipity is by no means trivial; for example, this na\"ive anti-telescoping fails to help show \eqref{finiteRRdiff} since, among numerous other terms, the coefficient of $q^8$ is $-1$ in the second ($i=2$) addend of the na\"ive anti-telescoping of \eqref{finiteRRdiff} for every $L>1$. A less na\"ive approach might sometimes be more beneficial, but for our purposes in this paper the na\"ive approach outlined above is sufficient.

Now clearly we could multiply every exponent in any inequality akin to \eqref{RRdiff}--\eqref{BGrdiff} by some common factor to obtain an inequality without $(1-q)$ as the leading factor in the denominator on the left; when looking at the partition-theoretic interpretation, this creates ``reducible'' examples (but examples nonetheless) where parts of size 1 are not needed to ``fill in the gaps''. In 2012, at the Ramanujan 125 Conference in Gainesville, Florida, Hamza Yesilyurt asked if the inclusion of the factor $(1-q)$ was necessary in all irreducible inequalities. We are pleased to answer in the negative, as stated in the following new theorem.
\begin{theorem}\label{xythm}
For any sextuple of positive integers $(L,m,x,y,r,R)$, 
\begin{equation*}
\frac{1}{(q^x,q^y,q^{rx+Ry};q^m)_L}\GEQ
\frac{1}{(q^{rx},q^{Ry},q^{x+y};q^m)_L}.
\end{equation*}
\end{theorem}
Clearly Theorem \ref{xythm} yields infinitely many irreducible examples. More astounding, however, is that the modulus $m$ can be \emph{arbitrary}. Even more amazing still is the relative ease with which the proof can be written using na\"ive anti-telescoping!

It is also possible, albeit more difficult, to use na\"ive anti-telescoping to yield the following new theorem.
\begin{theorem}\label{xyzthm}
For any octuple of positive integers $(L,m,x,y,z,r,R,\rho)$, 
\begin{equation*}
\frac{1}{(q^x,q^y,q^z,q^{rx+Ry+\rho z};q^m)_L}\GEQ
\frac{1}{(q^{rx},q^{Ry},q^{\rho z},q^{x+y+z};q^m)_L}.
\end{equation*}
\end{theorem}
The extra difficulty in proving Theorem \ref{xyzthm} comes from the fact that it seems to be impossible to re-write the addends in a natural way that makes it obvious that each addend only contributes nonnegative coefficients to the $q$-series. Consequently, en route to proving Theorem \ref{xyzthm}, we will require the following unobvious result, which is worthwhile in its own right and is not found anywhere else. (Most notably, we do not find anything of this form in \cite{IG},
which contains a compendium of rational functions with nonnegative coefficients.)

\begin{lemma}\label{mainlemma}
Let $r$ and $R$ be positive integers. Then the multivariate rational function 
\begin{equation*}
f(x,y,t) := \frac{(1-xy)(1-tx^r)(1-ty^R)+(1-t^2)(x-x^r)(y-y^R)}{(1-tx^r)(1-ty^R)(1-x)(1-y)(1-tx)(1-ty)},
\end{equation*}
with $|x|<1$, $|y|<1$, and $|t|<1$, has nonnegative coefficients when written as a power series centered at $(0,0,0)$.
\end{lemma}

In section \ref{sec:xythm proof}, we provide a proof of Theorem \ref{xythm} using a simple rational function identity together with na\"ive anti-telescoping, followed by a discussion of a partition theoretic interpretation of the difference 
\begin{equation*}
\frac{1}{(q^x,q^y,q^{rx+Ry};q^m)_L} - \frac{1}{(q^{rx},q^{Ry},q^{x+y};q^m)_L}.
\end{equation*}
In section \ref{sec:lemma proof} we give a proof of Lemma \ref{mainlemma}, which will be used in the proof of Theorem \ref{xyzthm} in section \ref{sec:xyzthm proof}. 
We then conclude in section \ref{sec:concluding remarks} with a brief discussion of a more general inequality.

\section{Proof of Theorem \ref{xythm}}\label{sec:xythm proof}

Let $P(i) := (q^x,q^y,q^{rx+Ry};q^m)_i$ and $Q(i) := (q^{rx},q^{Ry},q^{x+y};q^m)_i$.
We observe that since the identity
\begin{equation*}
\begin{split}
(1-t\alpha)(1-t\beta)(1-txy)-(1-tx)(1-ty)(1-t\alpha\beta)\\
=
t(x-\alpha)(1-\beta)(1-ty)+t(y-\beta)(1-t\alpha)(1-x)
\end{split}
\end{equation*}
is true, substituting $q^x$, $q^y$, $q^{rx}$, and $q^{Ry}$ for $x$, $y$, $\alpha$, and $\beta$, respectively, we can conclude that
\begin{equation}\label{TwoVarIDa}
(1-tq^{rx})(1-tq^{Ry})(1-tq^{x+y}) - (1-tq^x)(1-tq^y)(1-tq^{rx+Ry})
\end{equation}
and
\begin{equation}\label{TwoVarIDb}
tq^x(1-q^{(r-1)x})(1-q^{Ry})(1-tq^y) + tq^y(1-q^{(R-1)y})(1-tq^{rx})(1-q^x)
\end{equation}
are identically equal. 
Letting $t=q^{(i-1)m}$, we may use the equality of \eqref{TwoVarIDa} and \eqref{TwoVarIDb} to write
\begin{equation*}
\frac{Q(i-1)}{P(i)Q(L)}\left(\frac{Q(i)}{Q(i-1)}-\frac{P(i)}{P(i-1)}\right) 
=
V(i) + W(i),
\end{equation*}
where
\begin{equation*}
V(i) := \frac{q^{m(i-1)+y}(1-q^{(R-1)y})(1-q^x)(1-q^{m(i-1)+rx})}{P(i)\cdot Q(L)/Q(i-1)},
\end{equation*}
and
\begin{equation*}
W(i) := \frac{q^{m(i-1)+x}(1-q^{(r-1)x})(1-q^{Ry})(1-q^{m(i-1)+y})}{P(i)\cdot Q(L)/Q(i-1)}.
\end{equation*}
We note that since $(1-q^x)$ and $(1-q^y)$ are factors of the product $P(i)$ and since $(1-q^{m(i-1)+rx})$ is a factor of the product $Q(L)/Q(i-1)$, we have $V(i)\GEQ0$ for $1 \leq i \leq L$. To see that $W(i)\GEQ0$, we consider the following two cases.
\begin{enumerate}
\item Suppose $i=1$; then $(1-q^x)$ and $(1-q^y)=(1-q^{m(i-1)+y})$ are factors of $P(i)=P(1)$ and $(1-q^{Ry})$ is a factor of $Q(L)/Q(i-1) = Q(L)$. Thus, $W(1)\GEQ 0$.
\item Suppose $i>1$; then $(1-q^x)$, $(1-q^y)$, and $(1-q^{m(i-1)+y})$ are all independent factors of $P(i)$. Thus, $W(i)\GEQ0$.
\end{enumerate}
Finally, applying the anti-telescoping \eqref{antitelxy}, we have
\begin{equation}\label{V+W_sum}
\frac{1}{P(L)}-\frac{1}{Q(L)}
=
\sum_{i=1}^L \left( V(i) + W(i) \right),
\end{equation}
which then suffices to prove the theorem.
\qed

It would be nice to have a combinatorial proof of \eqref{V+W_sum}, but such has not been discovered by the time this paper was written. We note, however, that a partition interpretation of the right-hand side of \eqref{V+W_sum} is possible. Given a partition $\pi$, we let $p_j$ denote the part that is equal to $p+(j-1)m$, and we let $\num(p_j,\pi)$ represent the number of occurrences of the part $p_j$ in the partition $\pi$. Then, for a fixed $L$ we define
\begin{align*}
  \M(p,\pi) &:= \max \left(\{j : \num(p_j,\pi)>0\}\cup\{0\}\right)\\
\intertext{and}
  \m(p,\pi) &:= \min \left(\{j : \num(p_j,\pi)>0\}\cup\{L+1\}\right).
\end{align*}
We may consider $\sum_{i=1}^L V(i)$ and $\sum_{i=1}^L W(i)$, from \eqref{V+W_sum}, as two separate generating functions for partitions into parts congruent to (for $1 \leq i \leq L$) $x_i$, $y_i$, $(x+y)_i$, $(rx)_i$, $(ry)_i$, or $(rx+ry)_i$, subject to certain restrictions. (Note: in the cases where a particular part could arise in multiple ways, for example if $x_3=y_1$ or $rx=y$, then it would be necessary to treat the parts that arise in different ways as distinct, perhaps by assigning them unique colors based on what the base part is; since the base part is always one of $x$, $y$, $(x+y)$, $rx$, $Ry$, and $(rx+Ry)$, no more than six colors should be required.) 
We may take the restrictions as follows.
\begin{center}
\begin{tabular}{rlrl}
\multicolumn{2}{l}{Restrictions for $\sum_{i=1}^L V(i)$:}&
\multicolumn{2}{l}{Restrictions for $\sum_{i=1}^L W(i)$:}\\
V1: & $\M(y,\pi) \geq \max(\{1,\M(x,\pi)\})$ & W1: & $\M(x,\pi) > \M(y,\pi)$\\
V2: & $\M(y,\pi) \geq \M(rx+Ry,\pi)$         & W2: & $\M(x,\pi) \geq \M(rx+Ry,\pi)$\\
V3: & $\m(rx,\pi) > \M(y,\pi)$               & W3: & $\m(rx,\pi) \geq \M(x,\pi)$\\
V4: & $\m(Ry,\pi) \geq \M(y,\pi)$            & W4: & $\m(Ry,\pi) \geq \max(\{2,\M(x,\pi)\})$\\
V5: & $\m(x+y,\pi) \geq \M(y,\pi)$           & W5: & $\m(x+y,\pi) \geq \M(x,\pi)$\\
V6: & $\num(x_1,\pi) = 0$                    & W6: & $\num(x_1,\pi) < r-1$\\
V7: & $\num(y_1,\pi) < R-1$                  & W7: & $\num(y_1,\pi) < R$
\end{tabular}
\end{center}
Since the restrictions V1 and W1 are mutually exclusive, we may consider the right-hand side of \eqref{V+W_sum} as the generating function for partitions into parts congruent to (for $1 \leq i \leq L$) $x_i$, $y_i$, $(x+y)_i$, $(rx)_i$, $(Ry)_i$, or $(rx+Ry)_i$ such that the partition satisfies either V1--V7 or W1--W7.

\section{Proof of Lemma \ref{mainlemma}}\label{sec:lemma proof}

Let $[t^n]F(t)$ denote the coefficient of $t^n$ extracted from $F(t)$ (when written as a Maclaurin series). Direct calculations yield
\begin{equation}\label{eqone}
\begin{split}
[t^n]f(x,y,t) = 
& \phantom{{}+{}}  \frac{(1-xy)(x^{n+1}-y^{n+1})}{(1-x)(1-y)(x-y)}\\
& + \frac{(-x^{n+r}(1-x^2)+x^{nr+1}(1-x^{2r}))(y-y^R)}{(1-x)(1-y)(x-y)(x^r-y^R)}\\
& + \frac{(-y^{n+r}(1-y^2)+y^{nR+1}(1-y^{2R}))(x-x^r)}{(1-x)(1-y)(x-y)(x^r-y^R)}\\
& + \frac{(x^{(n-1)r}(1-x^{2r})-y^{n-1}(1-y^2))yx^r(x-x^r)(y-y^R)}{(1-x)(1-y)(x-y)(x^r-y^R)(x^r-y)}\\
& + \frac{(y^{(n-1)R}(1-y^{2R})-x^{n-1}(1-x^2))xy^R(x-x^r)(y-y^R)}{(1-x)(1-y)(x-y)(x^r-y^R)(y^R-x)}.
\end{split}
\end{equation} 
Claim:
\begin{equation}\label{eqtwo}
\begin{split}  
[t^n]f(x,y,t) = 
& \phantom{{}+{}}   \frac{x^n(1-y^{n+1})}{(1-y)(1-x)}
  + \frac{(y^{n+1}-y^{(n+1)R})(x^n-x^r)}{(1-y)(1-x)}\\
& + \frac{(y^n-y^{nR})(x^2-x^{2r})}{(1-y)(1-x)}
  + \frac{x(y^n-y^{(n+1)R})}{1-y}\\ 
& + \sum_{j=1}^{n-1} \frac{x^{(n-j)r}(y^j -y^{jR})(1-x^{2r})}{(1-y)(1-x)}\\
& + \sum_{j=0}^{(n-2-\delta(n))/2} \frac{x^{n-2j-1}y^{R(2j+1)}(1+x)}{1-y}\\
& + \sum_{j=1}^{(n-2+\delta(n))/2} \frac{x^{n-2j}y^{2jR}(1-y^{R(n+1-2j)})(1+x)}{1-y}\\
& + \frac{y^n}{1-y} 
  + \frac{\delta(n)xy^{(n+1)R}}{1-y},
\end{split}
\end{equation}
where $\delta(n)=0$ if $n$ is even and $\delta(n)=1$ if $n$ is odd.
To verify \eqref{eqtwo}, one first eliminates the sums in \eqref{eqtwo} to obtain
\begin{equation}\label{eqthree}
\begin{split}
[t^n]f(x,y,t) = 
& \phantom{{}+{}}  \frac{x^n(1-y^{n+1})}{(1-y)(1-x)}
  + \frac{(y^{n+1} - y^{(n+1)R})(x^n-x^r)}{(1-y)(1-x)}\\
& + \frac{(y^n-y^{nR})(x^2-x^{2r})}{(1-y)(1-x)}
  + \frac{x(y^n-y^{(n+1)R})}{1-y}\\
& + \frac{y^n}{1-y}
  + \frac{(1+x)xy^R(x^{n-1}-y^{(n-1)R})}{(1-y)(x-y^R)}\\
& + \frac{yx^r(x^{(n-1)r}-y^{n-1})(1-x^{2r})}{(1-y)(1-x)(x^r-y)}
  - \frac{y^{R(n+1)}(1+x)(x^2-x^n)}{(1-y)(1-x^2)}\\
& - \frac{y^Rx^r(x^{(n-1)r}-y^{R(n-1)})(1-x^{2r})}{(1-y)(1-x)(x^r-y^R)}.
\end{split} 
\end{equation}
Then, one can either verify by hand or use any number of symbolic manipulation programs to verify that the right-hand sides of \eqref{eqthree} and \eqref{eqone} are equal by simplifying their difference and getting 0. (The authors used Maple.)

We now observe that \eqref{eqtwo} implies that $[t^n]f(x,y,t)$ has nonnegative coefficients, provided $r \geq n$. Moreover, the only possible negative coefficients are
\[ 
[x^jy^k t^n]f(x,y,t) \quad \text{with} \quad 1< r < n \quad \text{and} \quad r \leq j < n < k < (n+1)R
\]
since all terms of \eqref{eqtwo} yield manifestly nonnegative coefficients except for the second term when $r<n$, where we have
\[
\dfrac{(y^{n+1}-y^{(n+1)R})(x^n-x^r)}{(1-y)(1-x)}
=
-(y^{n+1}+\cdots+y^{(n+1)R-1})(x^r+\cdots+x^{n-1}).
\]

Now suppose that the coefficient of $x^jy^kt^n$ in the power series for $f(x,y,t)$, centered at $(0,0,0)$, were negative; i.e.\ $[x^jy^kt^n]f(x,y,t)<0$.
Then, we must have both $1 < r < n$ and $r \leq j < n < k < R(n+1)$. 
Further, by the symmetry of $f(x,y,t)$ (with respect to the simultaneous swapping of $x$ and $r$ with $y$ and $R$, respectively) we would know that $[x^ky^jt^n]f(x,y,t)<0$ as well, and hence $1 < R < n$ and $R \leq k < n < j < r(n+1)$. 
But then we have a contradiction since we would have both $j<k$ and $k<j$. 
Thus, $[x^jy^kt^n]f(x,y,t)\geq0$, and the lemma is proved.
\qed

\section{Proof of Theorem \ref{xyzthm}}\label{sec:xyzthm proof}

Let $P(i) := (q^x,q^y,q^z,q^{rx+Ry+\rho z};q^m)_i$ and $Q(i) := (q^{rx},q^{Ry},q^{\rho z},q^{x+y+z};q^m)_i$.
Our goal will be to show that each addend in the sum on the right hand side of \eqref{antitelxyz} has nonnegative coefficients. We will do this by considering two cases based on the index of summation $i$ in \eqref{antitelxyz}: $i=1$ and $2 \leq i \leq L$.
First, though, we observe that 
\begin{equation}
(1-t\alpha)(1-t\beta)(1-t\gamma)(1-txyz)-(1-tx)(1-ty)(1-tz)(1-t\alpha\beta\gamma)
\end{equation}
is identically equal to
\begin{equation}
\begin{split}
&\phantom{{}+{}}\tfrac{1}{2}t(x-\alpha)\left[(1-t\beta)(1-t\gamma)(1-yz)+(1-ty)(1-tz)(1-\beta\gamma)\right]\\
&+\tfrac{1}{2}t(y-\beta)\left[(1-t\gamma)(1-t\alpha)(1-zx)+(1-tz)(1-tx)(1-\gamma\alpha)\right]\\
&+\tfrac{1}{2}t(z-\gamma)(1-tx)(1-ty)(1-\alpha\beta)\\
&+\tfrac{1}{2}t(z-\gamma)\left[(1-t\alpha)(1-t\beta)(1-xy)+(1-t^2)(x-\alpha)(y-\beta)\right].
\end{split}
\end{equation}
Substituting $q^x$, $q^y$, $q^z$, $q^{rx}$, $q^{Ry}$, and $q^{\rho z}$ for $x$, $y$, $z$, $\alpha$, $\beta$, and $\gamma$, respectively, we may then conclude that
\begin{equation}\label{main-pre-split}
\begin{split}
&\phantom{{}-{}}(1-tq^{rx})(1-tq^{Ry})(1-tq^{\rho z})(1-tq^{x+y+z})\\
&-(1-tq^x)(1-tq^y)(1-tq^z)(1-tq^{rx+Ry+\rho z})
\end{split}
\end{equation}
is identically equal to
\begin{equation}\label{main-split}
\begin{split}
&\phantom{{}+{}} \tfrac{1}{2}tq^x(1-q^{(r-1)x})\left[(1-tq^{Ry})(1-tq^{\rho z})(1-q^{y+z})\right.\\
 & \quad\quad\quad\quad\quad\quad\quad\quad\quad + \left.(1-tq^y)(1-tq^z)(1-q^{Ry+\rho z})\right]\\
&+ \tfrac{1}{2}tq^y(1-q^{(R-1)y})\left[(1-tq^{\rho z})(1-tq^{rx})(1-q^{z+x})\right.\\
 & \quad\quad\quad\quad\quad\quad\quad\quad\quad + \left.(1-tq^z)(1-tq^x)(1-q^{\rho z+rx})\right]\\
&+ \tfrac{1}{2}tq^z(1-q^{(\rho-1)z})(1-tq^x)(1-tq^y)(1-q^{rx+Ry})\\
&+ \tfrac{1}{2}tq^z(1-q^{(\rho-1)z})\left[(1-tq^{rx})(1-tq^{Ry})(1-q^{x+y})\right.\\
 & \quad\quad\quad\quad\quad\quad\quad\quad\quad + \left.(1-t^2)(q^x-q^{rx})(q^y-q^{Ry})\right].
\end{split}
\end{equation}

Let $t:=q^{(i-1)m}$.
Then, the numerator of the $i$th addend in \eqref{antitelxyz}, namely
\begin{equation*}
\frac{Q(i)}{Q(i-1)} - \frac{P(i)}{P(i-1)},  
\end{equation*}
is given precisely by \eqref{main-split}.
Now turning to the denominator of \eqref{antitelxyz}, we may write
\begin{equation*}
\frac{Q(L)}{Q(i-1)} = \frac{(q^{rx},q^{Ry},q^{\rho z},q^{x+y+z};q^m)_L}{(q^{rx},q^{Ry},q^{\rho z},q^{x+y+z};q^m)_{i-1}} = (tq^{rx},tq^{Ry},tq^{\rho z},tq^{x+y+z};q^m)_{L-i+1},
\end{equation*}
and so we have that 
\begin{equation}\label{Q(L)/Q(i-1)}
(1-tq^{rx})(1-tq^{Ry})(1-tq^{\rho z}) \text{ divides } \frac{Q(L)}{Q(i-1)} \text{ whenever } 1 \leq i \leq L.
\end{equation}
Similarly, from the definition of $P(i)$ we may deduce that
\begin{equation}\label{P(i),i=1}
(1-q^x)(1-q^y)(1-q^z) \text{ divides } P(1),
\end{equation}
and whenever $i>1$ that
\begin{equation}\label{P(i),i>1}
(1-q^x)(1-q^y)(1-q^z)(1-tq^x)(1-tq^y)(1-tq^z) \text{ divides } P(i).
\end{equation}

When $i=1$ we have $t=1$, and hence the numerator of the first addend in \eqref{antitelxyz} simplifies to
\begin{equation}\label{L=1split}
\begin{split}
Q(1)-P(1) = &\phantom{{}+{}} \tfrac{1}{2}q^x(1-q^{(r-1)x})\left[(1-q^{Ry})(1-q^{\rho z})(1-q^{y+z})\right.
\\ & \quad\quad\quad\quad\quad\quad\quad\quad\quad 
 + \left.(1-q^y)(1-q^z)(1-q^{Ry+\rho z})\right]\\
&+ \tfrac{1}{2}q^y(1-q^{(R-1)y})\left[(1-q^{\rho z})(1-q^{rx})(1-q^{z+x})\right.
\\ & \quad\quad\quad\quad\quad\quad\quad\quad\quad 
 + \left.(1-q^z)(1-q^x)(1-q^{\rho z + rx})\right]\\
&+ \tfrac{1}{2}q^z(1-q^{(\rho-1)z})\left[(1-q^x)(1-q^y)(1-q^{rx+Ry})\right.
\\ & \quad\quad\quad\quad\quad\quad\quad\quad\quad 
 + \left.(1-q^{rx})(1-q^{Ry})(1-q^{x+y})\right].
\end{split}
\end{equation}
Meanwhile, the denominator of the first addend in \eqref{antitelxyz} contains all of the factors indicated in \eqref{Q(L)/Q(i-1)}: $(1-q^{rx})$, $(1-q^{Ry})$, $(1-q^{\rho z})$. The denominator also contains all of the factors indicated by \eqref{P(i),i=1}: $(1-q^x)$, $(1-q^y)$, $(1-q^z)$. These factors, together with the ``trick'' of re-writing, for example,
\begin{equation}\label{trick}
(1-q^{x+y}) = (1-q^x)+q^x(1-q^y),
\end{equation}
is enough to see that the first addend in \eqref{antitelxyz} only has nonnegative coefficients.

When $2\leq i \leq L$, we have $t=q^{(i-1)m}\neq1$, and hence the numerator of the $i$th addend in \eqref{antitelxyz} is precisely \eqref{main-split}. From \eqref{Q(L)/Q(i-1)} and \eqref{P(i),i>1} we have the following factors in the denominator: $(1-tq^{rx})$, $(1-tq^{Ry})$, $(1-tq^{\rho z})$, $(1-q^x)$, $(1-q^y)$, $(1-q^z)$, $(1-tq^x)$, $(1-tq^y)$, $(1-tq^z)$. Again employing the ``trick'' \eqref{trick} as necessary, we can handle most of the $i$th addend similar to before, except for the last term of \eqref{main-split}, which contains the factor 
\begin{equation}\label{problemfactor}
\left[(1-q^{x+y})(1-tq^{rx})(1-tq^{Ry})+(1-t^2)(q^x-q^{rx})(q^y-q^{Ry})\right].
\end{equation}
This factor is potentially problematic due to the presence of the factor $(1-t^2)$ in the second term.

If we let $f$ be given as in Lemma \ref{mainlemma}, then \eqref{problemfactor} becomes
\begin{equation*}
f(q^x,q^y,t)(1-tq^{rx})(1-tq^{Ry})(1-q^x)(1-q^y)(1-tq^x)(1-tq^y).  
\end{equation*}
The last term of \eqref{main-split}, when divided by the nine factors listed above, then becomes
\begin{align*}
\frac{\tfrac{1}{2}tq^z(1-q^{(\rho-1)z})f(q^x,q^y,t)}{(1-tq^{\rho z})(1-q^z)(1-tq^z)},
\end{align*}
which, in light of Lemma \ref{mainlemma}, clearly now has no negative coefficients. Thus, having shown that all addends in \eqref{antitelxyz} admit only nonnegative coefficients, Theorem \ref{xyzthm} is proved.
\qed

\section{Concluding Remarks}\label{sec:concluding remarks}

It seems to always be possible to find a suitable ``splitting'' to handle the $L=1$ case, no matter how many variables are used. For example, if we increase from three to four main variables (adding in $w$ and $P w$ now), for $L=1$ we have
\begin{equation*}
\begin{split}
&\phantom{{}-{}}\frac{1}{(1-q^x)(1-q^y)(1-q^z)(1-q^w)(1-q^{rx+Ry+\rho z+P w)})}\\
& {}-{}\frac{1}{(1-q^{rx})(1-q^{Ry})(1-q^{\rho z})(1-q^{P w})(1-q^{x+y+z+w})}\\
&  = \frac{h(x,y,z,r,R,\rho)+h(x,y,w,r,R,P)+h(x,z,w,r,\rho,P)+h(y,z,w,R,\rho,P)}{(1-q^{x+y+z+w})(1-q^{rx+Ry+\rho z+P w})},
\end{split}
\end{equation*}
where
$h(x,y,z,r,R,\rho):=$
\begin{equation*}
\begin{split}
&\phantom{{}+{}} \tfrac{q^x(1-q^{(r-1)x})}{(1-q^x)(1-q^{rx})} \cdot \tfrac{q^y(1-q^{(R-1)y})}{(1-q^y)(1-q^{Ry}} \cdot \tfrac{q^z(1-q^{(\rho-1)z})}{(1-q^z)(1-q^{\rho z})}
+\tfrac{1}{2} \cdot \tfrac{q^x(1-q^{(r-1)x})}{(1-q^x)(1-q^{rx})} \cdot \tfrac{q^y(1-q^{(R-1)y})}{(1-q^y)(1-q^{Ry})}\\
&+\tfrac{1}{2} \cdot \tfrac{q^y(1-q^{(R-1)y})}{(1-q^y)(1-q^{Ry})} \cdot \tfrac{q^z(1-q^{(\rho-1)z})}{(1-q^z)(1-q^{\rho z})}
+\tfrac{1}{2} \cdot \tfrac{q^x(1-q^{(r-1)x})}{(1-q^x)(1-q^{rx})} \cdot \tfrac{q^z(1-q^{(\rho-1)z})}{(1-q^z)(1-q^{\rho z})}\\
&+\tfrac{1}{2} \cdot \tfrac{q^x(1-q^{(r-1)x})}{1-q^x} \cdot \tfrac{q^{Ry}}{(1-q^{rx})(1-q^{Ry})}
+\tfrac{1}{2} \cdot \tfrac{q^x(1-q^{(r-1)x})}{1-q^x} \cdot \tfrac{q^{\rho z}}{(1-q^{rx})(1-q^{\rho z})}\\
&+\tfrac{1}{2} \cdot \tfrac{q^y(1-q^{(R-1)y})}{(1-q^y)(1-q^{Ry})} \cdot \tfrac{q^{\rho z}}{1-q^{\rho z}}
+\tfrac{1}{2} \cdot \tfrac{q^y(1-q^{(R-1)y})}{(1-q^y)(1-q^{Ry})} \cdot \tfrac{q^{rx}}{1-q^{rx}}
+\tfrac{1}{2} \cdot \tfrac{q^z(1-q^{(\rho-1)z})}{(1-q^z)(1-q^{\rho z})} \cdot \tfrac{q^{Ry}}{1-q^{Ry}}\\
&+\tfrac{1}{2} \cdot \tfrac{q^z(1-q^{(\rho-1)z})}{(1-q^z)(1-q^{\rho z})} \cdot \tfrac{q^{rx}}{1-q^{rx}}
+\tfrac{1}{3} \cdot \tfrac{q^x(1-q^{(r-1)x})}{(1-q^x)(1-q^{rx})}
+\tfrac{1}{3} \cdot \tfrac{q^y(1-q^{(R-1)y})}{(1-q^y)(1-q^{Ry})}\\
&+\tfrac{1}{3} \cdot \tfrac{q^z(1-q^{(\rho-1)z})}{(1-q^z)(1-q^{\rho z})}
+\tfrac{q^x(1-q^{(r-1)x})}{(1-q^x)(1-q^{rx})} \cdot \tfrac{q^y(1-q^{(R-1)y})}{(1-q^y)(1-q^{Ry})} \cdot \tfrac{q^{\rho z}}{1-q^{\rho z}}\\
&+\tfrac{q^x(1-q^{(r-1)x})}{(1-q^x)(1-q^{rx})} \cdot \tfrac{q^z(1-q^{(\rho-1)z})}{(1-q^z)(1-q^{\rho z})} \cdot \tfrac{q^{Ry}}{1-q^{Ry}}
+\tfrac{q^y(1-q^{(R-1)y})}{(1-q^y)(1-q^{Ry})} \cdot \tfrac{q^z(1-q^{(\rho-1)z})}{(1-q^z)(1-q^{\rho z})} \cdot \tfrac{q^{rx}}{1-q^{rx}}\\
&+\tfrac{q^x(1-q^{(r-1)x})}{(1-q^x)(1-q^{rx})} \cdot \tfrac{q^{Ry}}{1-q^{Ry}} \cdot \tfrac{q^{\rho z}}{1-q^{\rho z}}
+\tfrac{q^y(1-q^{(R-1)y})}{(1-q^y)(1-q^{Ry})} \cdot \tfrac{q^{rx}}{1-q^{rx}} \cdot \tfrac{q^{\rho z}}{1-q^{\rho z}}\\
&+\tfrac{q^z(1-q^{(\rho-1)z})}{(1-q^z)(1-q^{\rho z})} \cdot \tfrac{q^{Ry}}{1-q^{Ry}} \cdot \tfrac{q^{rx}}{1-q^{rx}}
\end{split}
\end{equation*}
satisfies $h(x,y,z,r,R,\rho) \GEQ 0$.
Finding a suitable ``splitting'' with $t:=q^{(i-1)m}$ inserted into opportune locations, as we did in the proofs of the Theorems \ref{xythm} and \ref{xyzthm}, is a much more difficult task here. (We think of this as inserting the $t$'s since we wish to recover the $L=1$ case when we let $t=1$.) The authors of this manuscript do not currently possess such a ``splitting'' for this case. Nonetheless, the authors are fairly confident in the veracity of the following proposal.
\begin{proposal}
\label{proposal}
For any $(2n+2)$-tuple $\left(L,m,x_{(1)},\dots,x_{(n)},r_{(1)},\dots,r_{(n)}\right)$ of positive integers,
\begin{equation}\label{proposalineq}
\frac{1}{(q^{x_{(1)}},\dots,q^{x_{(n)}},q^{\Sigma};q^m)_L}
\GEQ
\frac{1}{(q^{r_{(1)}x_{(1)}},\dots,q^{r_{(n)}x_{(n)}},q^{\sigma};q^m)_L},
\end{equation}
where $\Sigma:=r_{(1)}x_{(1)}+\cdots+r_{(n)}x_{(n)}$ and $\sigma:=x_{(1)}+\cdots+x_{(n)}$.
\end{proposal}
We note that Proposal \ref{proposal} is true for $L=1$ since the right-hand side of \eqref{proposalineq} could be interpreted as the generating function for partitions into parts from the set $S:=\{x_{(1)},\dots,x_{(n)},\Sigma\}$ (parts with the same numeric value but distinct origins having different colors, thus ensuring $|S|=n+1$) such that for any such partition $\pi$, there is an integer $A$ with the property that
\begin{align*}
A &\equiv \num(x_{(1)},\pi) \pmod{r_{(1)}},\\
A &\equiv \num(x_{(2)},\pi) \pmod{r_{(2)}},\\
  &\ \ \vdots\\
A &\equiv \num(x_{(n)},\pi) \pmod{r_{(n)}},
\end{align*}
where $\num(p,\pi)$ is the number of occurrences of the part $p$ in the partition $\pi$. 
This set of partitions is a subset of the set of all partitions into parts from the set $S$, which is what the left-hand side of \eqref{proposalineq} would count. To see this clearly, we let $\pi'$ be a partition with parts from the set $S':=\{r_{(1)}x_{(1)},\dots,r_{(n)}x_{(n)},\sigma\}$ and let $\mu':=\min(\{\num(r_{(i)}x_{(i)},\pi') : 1 \leq i \leq n\})$. Then we can explicitly define an injection (for $L=1$) mapping $\pi' \mapsto \pi$ as follows:
\begin{align*}
\num(\Sigma,\pi) &:= \mu',\\
\num(x_{(i)},\pi) &:= r_{(i)}\cdot(\num(r_{(i)}x_{(i)},\pi')-\mu')+\num(\sigma,\pi').
\end{align*}
Clearly we can then choose $A=\num(\sigma,\pi')$. Now this mapping is invertible since if we let $\mu:=\min(\{\num(x_{(i)},\pi) : 1 \leq i \leq n\})$ we have
\begin{align*}
\num(\sigma,\pi') &= \mu, \\
\num(r_{(i)}x_{(i)},\pi') &= \frac{\num(x_{(i)},\pi)-\mu}{r_{(i)}}+\num(\Sigma,\pi).
\end{align*}
Thus, the proposal is proved for $L=1$.

Finally, we intend to explore possible connections with the recent work ``A $q$-rious positivity'' by S.\ Ole Warnaar and Wadim Zudilin (see \cite{WZ}). In particular, we are $q$uite $q$-rious as to how the validity of inequalities, like those in this paper, for $L=1$ might imply the validity for all positive $L$, a sentiment that seems echoed by the authors of \cite{WZ}.

\vspace{1em}

\noindent
{\it Acknowledgements}.
We would like to thank George Andrews and Wadim Zudilin for their interest and helpful discussions.

\providecommand{\bysame}{\leavevmode\hbox to3em{\hrulefill}\thinspace}

\end{document}